\documentclass[12pt]{article}
              \usepackage{amsfonts}
              \usepackage{amssymb}

\newtheorem{thm}{Theorem}

\newtheorem{lemma}[thm]{Lemma}

\def\be{\begin{eqnarray}}
\def\ee{\end{eqnarray}}
\def\bee{\begin{eqnarray*}}
\def\eee{\end{eqnarray*}}

\def\ts{\textstyle}
\def\bra{\langle}
\def\ket{\rangle}

\def\rt2{\ts \frac{1}{\sqrt{2}} }

\def\ot{\otimes}
   \def\tr{\hbox{Tr} \,}

\def\bmat{\begin{pmatrix}}
\def\emat{\end{pmatrix}}

\def\qed{\vskip0.1in \noindent \fbox{\bf QED} \vskip0.2in}

\title{Comparison of matrix norms on bipartite spaces}

    \author{Christopher King and Nilufer Koldan\\
    \\
{\small      Department of Mathematics} \\
{\small      Northeastern University} \\
{\small      Boston MA 02115}}

\begin{document}

     \maketitle

     \begin{abstract}
Two non-commutative versions
of the classical $L^q(L^p)$ norm on the product matrix algebras
${\cal M}_n \otimes {\cal M}_m$ are compared. The first norm was defined recently by Carlen and Lieb,
as a byproduct of their analysis of certain convex functions on matrix spaces.
The second norm was defined by Pisier and others using results from the theory of operator spaces.
It is shown that the second norm is upper bounded by a constant multiple of the first
for all $1 \le p \le2$, $q \ge 1$. In one case
($2 = p < q$) it is also shown that there is no such lower bound, and hence that the norms are
inequivalent. It is conjectured that the norms are inequivalent in all cases.
       \end{abstract}
%
%

\section{Introduction}
Let ${\cal M}_n$ denote the algebra of $n \times n$ complex-valued matrices.
The Schatten norm on ${\cal M}_n$ provides a non-commutative version
of the classical $L^p$ norm. It is defined for $p \ge 1$ by
\be
|| A ||_p = \Big( \tr |A|^p \Big)^{1/p}
\ee
Many of the standard properties of the classical norm extend to the
Schatten norm,
including monotonicity, convexity, H\"older's inequality and duality.

In this paper we compare two non-commutative versions
of the classical $L^q(L^p)$ norm. The first version was introduced by Carlen and Lieb
in their recent paper \cite{CL2}, where it arose out of ideas connected with the central
theme of strong subadditivity of entropy. The second version arises from the work
of Pisier and others on operator spaces \cite{Pis1, Pis2}. Norms of this second type were used in the paper
\cite{DJKR} to prove results about completely positive maps on matrix
algebras. Thus both norms are connected to important ideas in quantum information theory,
and this motivates our study of the similarities and differences between them.

Recall that in the classical (commutative) case
a function $f : \mathbb{R}^2 \rightarrow \mathbb{R}$,
$(x,y)  \mapsto f(x,y)$, may also be regarded as a map
$f \,:\, \mathbb{R} \rightarrow L^p(\mathbb{R})$ , $x  \mapsto f(x,
\cdot) \equiv f_x(\cdot)$. The $L^q$ norm of this
map is
\be\label{class-norm}
\bigg( \int \| f_x \|_{p}^q d x \bigg)^{1/q} = \bigg( \int \Big( \int
|f(x,y)|^{p} d y \Big)^{q/p}  d x  \bigg)^{1/q}
\ee
and this defines the $L^q(L^p)$ norm of $f$.

The right side of (\ref{class-norm}) suggests a possible non-commutative version of this norm
for the bipartite space ${\cal M}_{n m} = {\cal M}_n \otimes {\cal M}_m$,
namely $\bigg( \tr_1 ( \tr_2 |Y|^p)^{q/p} \bigg)^{1/q}$, where $\tr_1$ is the partial trace over the
first or `outside' space ${\cal M}_n$, and $\tr_2$ is the
partial trace over the second or `inside' space ${\cal M}_m$.  While this specific expression does not
lead to a norm,
Carlen and Lieb  have recently constructed a norm based on this idea.
For positive semidefinite $Y$ they define
\be\label{def:Psi}
\Psi_{p,q}(Y) = \bigg( \tr_1 ( \tr_2 Y^p)^{q/p} \bigg)^{1/q}
\ee
In Theorem 1.1 in the paper \cite{CL2} it is proved that
for all $1 \le p \le 2$ and $q \ge 1$ the function $\Psi_{p,q}$ is
convex on the set of positive semidefinite matrices in ${\cal M}_n \otimes {\cal M}_m$.

Building on the convexity properties of (\ref{def:Psi}), Carlen
and Lieb defined a new norm on ${\cal M}_n \otimes {\cal M}_m$ in
the following way. First define for Hermitian matrices $X$ the
quantity
\be
||| X |||_{p,q} = \inf_{A,B} \Big\{ \Psi_{p,q}(A) +
\Psi_{p,q}(B) \,\Big|\, X + A = B, \,\, A \ge 0, B \ge 0 \Big\}
\ee
Then the Carlen-Lieb norm is defined for $1 \le p
\le 2$, $q \ge 1$ as
\be\label{CL}
\| Y \|_{CL} = \frac{1}{2} \,\,
\bigg| \bigg|\bigg| \bmat {0 & Y \cr Y^* & 0 } \emat
\bigg|\bigg|\bigg|_{p,q}
\ee
On the right side of (\ref{CL}) the
function $\Psi_{p,q}$ is applied to a matrix in ${\cal M}_{2nm}$.
This space must be split as a bipartite space in order to apply
the definition. Carlen and Lieb choose the splitting ${\cal
M}_{2nm}  = {\cal M}_{n} \ot ({\cal M}_{m} \ot  {\cal M}_{2} )$.
In other words, in the definition (\ref{def:Psi}) the inside space
${\cal M}_m$  is replaced by ${\cal M}_{m} \ot  {\cal M}_{2}$, and
the outside space is still ${\cal M}_{n}$.

A different approach to the question of defining norms for these
bipartite matrix spaces arises from the work of Effros and Ruan
\cite{ER1, ER2}, Pisier \cite{Pis1, Pis2} and Junge \cite{Ju1} on operator spaces.
Several alternative formulations of such norms were used in the paper \cite{DJKR}
to analyze CB norms of completely positive maps.
In this paper we will present and analyze these Pisier-type norms using matrix analytic methods,
without relying on results from the general theory of operator spaces.
Define the value $r$ by
\be\label{def:r}
\frac{1}{r} = \max \Big\{\frac{1}{p}, \frac{1}{q} \Big\} - \min \Big\{\frac{1}{p}, \frac{1}{q} \Big\}
\ee
The expression for the norm differs for the cases $p < q$ and $p > q$ (the subscript `NC'
stands for `non-commutative'). The first expression (\ref{Pis1}) was used  in \cite{DJKR},
while the second expression (\ref{Pis2}) is a modified version of one used in \cite{DJKR}
(the expression used in \cite{DJKR} had only one term inside the infimum).

\vskip0.3in

\noindent \fbox{$1 \le p \le q \le \infty$}
\vskip-0.2in
\be\label{Pis1}
\| Y \|_{NC}  = \sup_{A,B \in {\cal M}_n} \bigg\{\, \frac{\| (A \otimes I_m)
Y ( B \otimes I_m) \|_{p}}{\| A \|_{2r} \, \| B \|_{2r}}  \bigg\}
\ee

\vskip0.1in

\noindent \fbox{$1 \le q \le p \le \infty$}
\vskip-0.2in

\be\label{Pis2}
\| Y \|_{NC}    =  \inf_{\stackrel{\scriptstyle A_i,B_i \in {\cal M}_n}{Z_i  \in {\cal M}_{nm}}} \Big\{\, \sum_i \| A_i \|_{2r} \, \| B_i \|_{2r}
\| Z_i \|_{p}
 : Y = \sum_i (A_i \ot I_m) Z_i (B_i \ot I_m) \Big\}
\ee

\bigskip
\noindent \underline{\em Remark  1:} without loss of generality the $\sup$ on the right side of (\ref{Pis1}) may be restricted to
the set of positive semidefinite matrices $A,B \ge 0$. This may be seen by writing $A = U (A^* A)^{1/2}$ and $B = (B B^*)^{1/2} V$ where
$U,V$ are unitaries, and observing that $\| U C \|_{t} = \| C U \|_t = \| C \|_t$ for all $t$, $C$ and unitary $U$.

\bigskip
\noindent \underline{\em Remark  2:} for positive semidefinite matrices $Y \ge 0$,
\be
\| (A \otimes I_m) Y ( B \otimes I_m) \|_{p} &= & \| (A \otimes I_m)
Y^{1/2} Y^{1/2} ( B \otimes I_m) \|_{p} \nonumber \\
& \le & \| (A \otimes I_m) Y ( A^{*} \otimes I_m) \|_{p}^{1/2} \,
\| (B^* \otimes I_m) Y ( B \otimes I_m) \|_{p}^{1/2} \nonumber
\ee
Thus the supremum on the right side of (\ref{Pis1}) may be restricted to $A=B \ge 0$.
Furthermore, letting
\be\label{def:C}
C = A^{2r} \,\, \bigg(\| A \|_{2r}\bigg)^{-2r}
\ee
it follows that $\tr C = 1$, and therefore for $p \le q$
\be
\| Y \|_{NC}  & = & \sup_{A \ge 0 \in {\cal M}_n} \bigg\{\, \frac{\| (A \otimes I_m)
Y ( A \otimes I_m) \|_{p}}{\| A \|_{2r}^2}  \bigg\} \\
& = & \sup_{C \ge 0, \, {\rm Tr} C = 1} \bigg\{\, \| (C^{1/2r} \otimes I_m)
Y ( C^{1/2r} \otimes I_m) \|_{p} \bigg\}
\ee

\bigskip
\noindent \underline{\em Remark  3:} in some cases we will need to
refer to the values $p,q$ in the norm; we do this by writing $\| \cdot \|_{NC:p,q}$.

\bigskip
Our first result establishes basic properties of $\| Y \|_{NC}$. These properties help to motivate the definitions (\ref{Pis1}) and (\ref{Pis2}).

\begin{lemma}\label{lem1}
Assume $p,q \ge 1$.
\par\noindent{\bf a)} [Triangle inequality] For any $Y,W \in {\cal M}_{nm}$
\be
\| Y + W \|_{NC} \le \| Y \|_{NC} + \| W \|_{NC}
\ee
\par\noindent{\bf b)} [H\"older's inequality] Define the usual conjugate values for $p,q$:
\bee
\frac{1}{p'} = 1 - \frac{1}{p}, \qquad \frac{1}{q'} = 1 - \frac{1}{q}
\eee
Then for all $Y,W \in {\cal M}_{nm}$,
\be
| \tr (Y W) | \le \| Y \|_{NC:p,q} \,\, \| W \|_{NC:p',q'}
\ee
\par\noindent{\bf c)} [Duality] For all $Y \in {\cal M}_{nm}$
\be
\| Y \|_{NC:p,q} = \sup_{W \in {\cal M}_{nm}} \, \{ | \tr (YW) | \, : \, \| W \|_{NC:p',q'} \le 1 \}
\ee
\par\noindent{\bf d)} [Product form for product matrices] For any $Y_1 \in {\cal M}_n$, $Y_2 \in {\cal M}_m$
\be
\| Y_1 \ot Y_2 \|_{NC} = \| Y_1 \|_{q} \, \| Y_2 \|_{p}
\ee
\par\noindent{\bf e)} [Reduction to Schatten norm at $p=q$] For any $Y \in {\cal M}_{nm}$, and for $p=q$,
\be
\| Y \|_{NC} = \| Y \|_{p}
\ee
\par\noindent{\bf f)} [Reduction to classical norm on diagonal matrices]
Choose orthonormal bases $\{e_i\}$ in $\mathbb{C}^n$ and $\{f_j\}$ in $\mathbb{C}^m$.
Then for any matrix $Y \in {\cal M}_{nm}$ which is diagonal in the basis $\{e_i \ot f_j \}$,
\be
\| Y \|_{NC} = \bigg( \tr_1 ( \tr_2 |Y|^p)^{q/p} \bigg)^{1/q}
\ee
\end{lemma}

\bigskip
In order for this paper to be self-contained, we include the proof of
Lemma \ref{lem1} in the Appendix.

\medskip
The main purpose of this paper is to compare the norms $\| \cdot \|_{CL}$ and $\| \cdot \|_{NC}$,
in particular to investigate whether the norms are equivalent.
The next result provides a  bound for $\| Y \|_{NC}$ in terms of
$\| Y \|_{CL}$.

\begin{thm}\label{thm1}
Assume $1 \le p \le 2$ and $q \ge 1$.
\par\noindent{\em a)}
Let $Y \ge 0$ be positive semidefinite, then
\be\label{thm1:bd1}
\| Y \|_{NC} \le \Psi_{p,q}(Y)
\ee
\medskip
\par\noindent{\em b)} For all matrices $Y$ in ${\cal M}_n \otimes {\cal M}_m$,
\be\label{thm1:bd2}
\| Y \|_{NC} \le 2^{3 - 1/p} \,\, \| Y \|_{CL}
\ee
\end{thm}

\medskip
\noindent \underline{\em Remark 4:} the bound (\ref{thm1:bd1}) is sharp, since both sides agree
on product matrices $Y = Y_1 \ot Y_2$.

\medskip
Recall that a function $f : {\cal M}_n \rightarrow \mathbb{R}$, $A \mapsto f(A)$ is monotone if
$0 \le A \le B$ implies $f(A) \le f(B)$.
In the paper \cite{CL2}, Remark 1.5, the authors point out that the function $A \mapsto \Psi_{p,q}(A)$ is not monotone,
but then state that the CL-norm is monotone. In fact as the following Lemma shows the CL-norm is
non-monotone in some cases.

\begin{lemma}\label{lem:non-mon}
For all $1 \le p \le 2$ and $p \le q$ the function $A \mapsto \| A \|_{CL}$ is not monotone.
\end{lemma}

\medskip
\noindent \underline{\em Remark 5:} 
In the case $p \le q$ it is clear from (\ref{Pis1}) that the NC-norm is monotone.
Thus Lemma \ref{lem:non-mon} implies that the CL- and NC-norms are different
in this case. The following Lemma addresses the question of whether the norms
are equivalent.

\begin{thm}\label{thm2}
Assume $p < q$.
\par\noindent{\em a)} Let $1 \le p \le 2$.
For all $k \ge 1$ there exist integers $\{n_k, m_k\}$ and positive semidefinite
matrices $0 \le Y^{(k)} \in {\cal M}_{n_k} \otimes {\cal M}_{m_k}$, such that
\be\label{thm2:ineq1}
\frac{\Psi_{p,q}(Y^{(k)})}{\| Y^{(k)} \|_{NC}} \rightarrow \infty \quad
{\rm as} \,\,\,\, k \rightarrow \infty
\ee
\medskip
\par\noindent{\em b)} Let $p = 2$. Then for all $k \ge 1$ there exist integers $\{n_k, m_k\}$ and
matrices $Y^{(k)} \in {\cal M}_{n_k} \otimes {\cal M}_{m_k}$, such that
\be\label{thm2:ineq2}
\frac{\| Y^{(k)} \|_{CL}}{\| Y^{(k)} \|_{NC}} \rightarrow \infty \quad
{\rm as} \,\,\,\, k \rightarrow \infty
\ee
\end{thm}

\bigskip
Theorem \ref{thm2} (b) implies that the
norms $\| Y \|_{NC}$ and $\| Y \|_{CL}$ are not equivalent when $p=2$ and $q > 2$, in the sense that
there do not exist non-zero finite constants $c_1$ and $c_2$ such that
$\| Y \|_{NC} \le c_1 \, \| Y \|_{CL} \le c_2 \, \| Y \|_{NC}$ for all matrices
$Y$. We conjecture that the norms are not equivalent for all
$1 < p \le 2$ and all $q \ge 1$.

\bigskip
The paper is organized as follows. In Section 2 we derive an alternative expression for the Carlen-Lieb
norm for the case where the matrix is Hermitian. In Section 3 we use this expression to first prove Theorem \ref{thm1}
for positive semidefinite matrices, and then extend the result to general matrices. The main technical tool
in the proof is an application of the Lieb-Thirring inequality. Section 4 presents
a construction of the counterexamples which prove Theorem \ref{thm2}, and the Appendix contains the 
proofs of Lemma \ref{lem1} and Lemma \ref{lem:non-mon}.

\section{Representation for CL-norm of Hermitian matrix}
In this section we derive a simplified representation for the Carlen-Lieb norm
in the case where the matrix is Hermitian. We assume throughout this section
that $1 \le p \le 2$ and $q \ge 1$.

\begin{lemma}\label{lem2}
For a Hermitian matrix $Y = Y^* \in {\cal M}_{mn}$,
\be\label{eqn:lem1}
\| Y \|_{CL} = \inf_{A \ge 0, \, Y+A \ge 0} \,\, \Psi_{p,q} \bmat {Y + A
& 0 \cr 0 & A} \emat
\ee
\end{lemma}

\par\noindent{\em Proof:}
recall the Pauli matrix $\sigma_x = \bmat {0 & 1 \cr 1 & 0} \emat$, then
\be
\bmat {0 & Y \cr Y & 0} \emat = Y \ot \sigma_x
\ee
Define
\be
{\cal G}(Y) = \{ A \in {\cal M}_{mn} \ot \mathbb{C}^2 \,|\, A \ge 0,
\,\, Y \ot \sigma_x + A \ge 0 \}
\ee
Then
\be\label{eqn1}
\| Y \|_{CL} = \frac{1}{2} \,\, \inf_{A \in {\cal G}(Y)} \,\, \bigg(
\Psi_{p,q}(Y \ot \sigma_x + A) + \Psi_{p,q}(A) \bigg)
\ee

In general $A \in {\cal G}(Y)$ has the form
$A = A_1 \ot I_2 + A_2 \ot \sigma_x + A_3 \ot \sigma_y + A_4 \ot
\sigma_z$, where $\sigma_y, \sigma_z$ are the
other Pauli matrices.
The function $\Psi_{p,q}$ is invariant under unitary transformations
on each of the spaces ${\cal M}_m$ and ${\cal M}_n$.
Hence $\Psi_{p,q}(A)$ and $\Psi_{p,q}(Y \ot \sigma_x + A)$ are
unchanged if $A$ is replaced
by $A' = (I \ot \sigma_x) A (I \ot \sigma_x) = A_1 \ot I_2 + A_2 \ot
\sigma_x - A_3 \ot \sigma_y - A_4 \ot \sigma_z$.
The function $\Psi_{p,q}$ is also convex, and hence the expression
$\Psi_{p,q}(Y \ot \sigma_x + A) + \Psi_{p,q}(A)$ on the right side of
(\ref{eqn1}) can only be lowered by
replacing $A$ by $(A + A')/2 = A_1 \ot I_2 + A_2 \ot \sigma_x$.
Therefore the infimum on the right
side of (\ref{eqn1}) may be restricted to matrices of the form
$A = A_1 \ot I + A_2 \ot \sigma_x$, where $A_1,A_2$ satisfy
\be
A_1 \ge | A_2 |, \quad
A_1 \ge | Y + A_2 |
\ee
The matrix $A$ is unitarily equivalent to $A_1 \ot I_2 + A_2 \ot
\sigma_z$, and similarly $Y \ot \sigma_x + A$ is unitarily equivalent to
$A_1 \ot I_2 + (Y + A_2) \ot  \sigma_z$.
Let $C = A_1 + A_2 \ge 0$ and $D = A_1 - A_2 \ge 0$.
Then (\ref{eqn1}) can be written as
\be
\| Y \|_{CL} = \frac{1}{2} \,  \inf_{\{C,D,C+Y,D-Y
\ge 0\}} \bigg(
\Psi \bmat { C + Y & 0 \cr 0 & D - Y} \emat +  \Psi \bmat { C & 0 \cr 0 &
D} \emat \bigg)
\ee

Alternatively, since $D = Y + E$ with $E \ge 0$ this can be written as
\be\label{eqn2}
\| Y \|_{CL}  = \frac{1}{2} \,  \inf_{\{C,E,C+Y,E+Y
\ge 0\}} \bigg(
\Psi \bmat { C + Y & 0 \cr 0 & E} \emat +  \Psi \bmat { C & 0 \cr 0 &
E + Y} \emat \bigg)
\ee
Note that
\be
\Psi \bmat { C & 0 \cr 0 &
E + Y} \emat = \Psi \bmat { E + Y & 0 \cr 0 &
C} \emat
\ee
and so (\ref{eqn2})  may be written as
\be\label{eqn3}
\| Y \|_{CL} = \frac{1}{2} \,  \inf_{\{C,E \ge 0, \,C+Y \ge 0, \,E+Y
\ge 0\}} \bigg(
\Psi \bmat { C + Y & 0 \cr 0 & E} \emat +  \Psi \bmat { E + Y & 0 \cr 0 &
C} \emat \bigg)
\ee
Convexity of $\Psi_{p,q}$ implies that
\be
\frac{1}{2} \, \bigg(
\Psi \bmat { C + Y & 0 \cr 0 & E} \emat +  \Psi \bmat { E + Y & 0 \cr 0 &
C} \emat  \bigg) \ge \Psi \bmat { F + Y & 0 \cr 0 & F} \emat
\ee
where $F = (C+E)/2$. Hence the infimum on the right side of (\ref{eqn3}) may be
restricted to $C=E$ and this gives (\ref{eqn:lem1}).

\section{Proof of Theorem \ref{thm1}}
\subsection{The bound for $Y \ge 0$}
\par\noindent \fbox{The case $p < q$}
Note that
\be\label{psi0.5}
\Psi_{p,q}(Y)^p = \bigg( \tr_1 ( \tr_2 Y^p)^{q/p} \bigg)^{p/q}
\ee
Define $t$ by
\be
\frac{1}{t} + \frac{p}{q} = 1
\ee
then by duality
\be\label{psi1}
\Psi_{p,q}(Y)^p = \sup_B \tr_1 \bigg[ B^{1/t} \, \tr_2 (Y^p) \bigg]
\ee
where the supremum runs over positive matrices $B \ge 0$ satisfying
$\tr B = 1$.
Let $r = pt$ then
\be
\frac{1}{r} + \frac{1}{q} = \frac{1}{p}
\ee
and  (\ref{psi1}) can be written as
\be\label{psi2}
\Psi_{p,q}(Y)^p = \sup_B \tr_1 \bigg[ B^{p/r} \, \tr_2 (Y^p) \bigg]
= \sup_B \tr \bigg[ (B^{p/r} \ot I_m) Y^p \bigg]
\ee

Turning now to the NC norm, by Remark 2
this may be rewritten as a supremum over positive matrices $B \ge 0$ with $\tr B =1$:
\be
\| Y \|_{NC}
= \sup_{B \ge 0, {\rm Tr} B = 1} \, \| (B^{1/2r} \otimes I_m) Y ( B^{1/2r} \otimes I_m) \|_{p}
\ee
Hence
\be
\| Y \|_{NC}^{p} = \sup_{B \ge 0, {\rm Tr} B = 1} \, \tr \bigg[ (B^{1/2r} \otimes I_m) Y
( B^{1/2r} \otimes I_m) \bigg]^p
\ee
The Lieb-Thirring inequality \cite{LiebTh} implies that for $Y,B \ge 0$ and $p \ge 1$,
\be
\tr \bigg[ (B^{1/2r} \otimes I_m) Y ( B^{1/2r} \otimes I_m) \bigg]^p \le
\tr \bigg[ (B^{p/2r} \otimes I_m) Y^p ( B^{p/2r} \otimes I_m) \bigg]
\ee
Comparing with (\ref{psi2}) yields the inequality.

\medskip
\par\noindent \fbox{The case $p > q$}
Given $0 < t < 1$ define the conjugate value $t'$ by
\be
\frac{1}{t'} = \frac{1}{t} - 1
\ee
Then for any non-negative sequence $a_1,\dots,a_n$ it is easy to check that
\be\label{dual:t<1}
\bigg( \sum_{i=1}^n a_i^t \bigg)^{1/t} = \inf \, \bigg\{ \sum_{i=1}^n a_i b_i \,:\,
b_i \ge 0, \, \sum_{i=1}^n b_i^{-t'} = 1 \bigg\}
\ee
Given a positive semidefinite matrix $A \ge 0$, $A \in {\cal M}_n$, define its positive commutant
\be
{\rm Comm}_{+}[A] = \{ B \in {\cal M}_n \,:\, B \ge 0, \, AB = BA \}
\ee
Then (\ref{dual:t<1}) implies that for $A \ge 0$
\bee
\| A \|_t & = & \inf \, \{ \tr (A B) \, :\, \tr B^{-t'} =1, \, B \in {\rm Comm}_{+}[A] \} \\
& = & \inf \, \{ \tr (A C^{-1/t'}) \, :\, \tr C =1, \, C \in {\rm Comm}_{+}[A] \}
\eee
Applying this to (\ref{def:Psi}) with $t = q/p$ and $A = \tr_2 (Y^p)$ gives
\be
\Psi_{p,q}(Y)^p  &= &  \inf_{{\rm Tr} \,C =1,\, C \in {\rm Comm}_{+}[A]} \, \tr_1 [C^{-p/r} \, \tr_2(Y^p) ] \nonumber \\
&=&  \inf_{{\rm Tr}\, C =1,\, C \in {\rm Comm}_{+}[A]} \,\tr [(C^{-p/2r} \ot I_m) Y^p (C^{-p/2r} \ot I_m) ]  \nonumber \\
& \ge & \inf_{{\rm Tr}\, C =1,\, C \ge 0} \,\tr [(C^{-p/2r} \ot I_m) Y^p (C^{-p/2r} \ot I_m) ] \label{psi2.5}
\ee

We obtain an upper bound for $\| Y \|_{NC}$ by
restricting the infimum on the right side of (\ref{Pis2}) to a single term $Y = (A \ot I_m) Z (B \ot I_m)$
with $A = B > 0$, which leads to the bound
\be\label{Pis2.5}
\| Y \|_{NC}^{p} \le \inf_{{\rm Tr}\, B =1,\, B \ge 0} \, \tr \bigg[ (B^{-1/2r} \otimes I_m) Y
( B^{-1/2r} \otimes I_m) \bigg]^p
\ee
Since $p > 1$ the Lieb-Thirring bound again implies that (\ref{psi2.5}) upper bounds
(\ref{Pis2.5}) and thus the result follows.

\subsection{The bound for general $Y$}
\par\noindent \fbox{The case $p < q$}
First we establish the bound for Hermitian matrices. Suppose $Y=Y^*$ then
\be\label{first:p<q}
\| Y \|_{NC} = \sup_{A,B \ge 0} \, \frac{\| (A \otimes I_m)
Y ( B \otimes I_m) \|_{p}}{\| A \|_{2r} \, \| B \|_{2r}}
\ee
Write $Y = Z - W$ where both $Z,W$ are positive, then
\be
&&\hskip-0.5in  \| (A \otimes I_m)
Y ( B \otimes I_m) \|_{p} \nonumber \\
&&  = \hskip0.15in \| (A \otimes I_m) Z ( B \otimes I_m) -  (A \otimes I_m) W ( B \otimes I_m) \|_{p} \nonumber \\
&& \le \hskip0.15in
\| (A \otimes I_m) Z ( B \otimes I_m)\|_{p}  +  \| (A \otimes I_m) W ( B \otimes I_m) \|_{p} \nonumber \\
&& \le \hskip0.15in
\| (A \otimes I_m) Z^{1/2} \|_{2p} \, \|  ( B \otimes I_m) Z^{1/2} \|_{2p} \nonumber \\
&& \hskip0.3in  +\hskip0.15in \| (A \otimes I_m) W^{1/2} \|_{2p} \, \| ( B \otimes I_m) W^{1/2}  \|_{2p} \nonumber \\
&& \le  \hskip0.15in \bigg( \| (A \otimes I_m) Z^{1/2} \|_{2p} +  \| (A \otimes I_m) W^{1/2} \|_{2p} \bigg) \nonumber \\
&& \hskip0.5in \,\, \times \,\, \bigg( \| ( B \otimes I_m) Z^{1/2}  \|_{2p} +  \|  ( B \otimes I_m) W^{1/2} \|_{2p} \bigg)
\ee
Taking the supremum over $A,B$ gives the bound
\be
\| Y \|_{NC} \le \sup_{A} \bigg( \| (A \otimes I_m) Z^{1/2} \|_{2p} +  \| (A \otimes I_m) W^{1/2} \|_{2p} \bigg)^2 \, (\| A \|_{2r})^{-2}
\ee
As noted in Remark 1 we can assume that $A \ge 0$. Let $C$ be defined as in (\ref{def:C}), then
\be
\| Y \|_{NC} & \le & \sup_{C} \bigg( \| (C^{1/2r} \otimes I_m) Z^{1/2} \|_{2p} +  \| (C^{1/2r} \otimes I_m) W^{1/2} \|_{2p} \bigg)^2 \nonumber \\
& \le &  2 \, \sup_{C} \bigg(  \| (C^{1/2r} \otimes I_m) Z^{1/2} \|_{2p}^2 + \| (C^{1/2r} \otimes I_m) W^{1/2} \|_{2p}^2 \bigg)
\ee
where the supremum is taken over positive matrices with $\tr C =1$. Note that
\be
\| (C^{1/2r} \otimes I_m) Z^{1/2} \|_{2p}^2 = \| (C^{1/2r} \otimes I_m) Z (C^{1/2r} \otimes I_m)  \|_{p}
\ee
and further that for all $x,y \ge 0$ we have the inequality
\be
x^{1/p} + y^{1/p} \le 2^{1 - 1/p} \, (x + y)^{1/p}
\ee
Hence
\be
\| Y \|_{NC}^p & \le & 2^{2p - 1} \, \sup_{C} \bigg( \tr \Big[(C^{1/2r} \otimes I_m) Z (C^{1/2r} \otimes I_m)\Big]^p \\
&& \hskip 1in + \tr \Big[(C^{1/2r} \otimes I_m) W (C^{1/2r} \otimes I_m)\Big]^p \bigg) \nonumber \\
& = & 2^{2p - 1} \, \sup_{C} \tr
\Big[(C^{1/2r} \otimes I_{2m}) \bmat {Z & 0 \cr 0 & W} \emat (C^{1/2r} \otimes I_{2m})\Big]^p \label{ineq4} \\
& = & 2^{2p - 1} \,\Big|\Big| \bmat {Z & 0 \cr 0 & W} \emat \Big|\Big|_{NC}^p \nonumber \\
& \le & 2^{2p - 1} \,\Psi_{p,q} \bmat {Z & 0 \cr 0 & W} \emat^p
\ee
where the norm on the right side is computed for the decomposition ${\cal M}_{2nm} = {\cal M}_{n} \otimes {\cal M}_{2m}$,
and where we used Theorem \ref{thm1} (a) at the last step. Since this holds for every pair of positive matrices
$Z,W$ satisfying $Y = Z - W$ we have by Lemma \ref{lem2}
\be\label{ineq6a}
\| Y \|_{NC} \le 2^{2 - 1/p} \, \inf_{W \ge 0, \, Y+W \ge 0} \, \Psi_{p,q} \bmat {Y+W & 0 \cr 0 & W} \emat
= 2^{2 - 1/p} \, \|Y\|_{CL}
\ee

This establishes the bound for Hermitian matrices. Now consider a general matrix $Y$, and write
$Y = Y_1 + i Y_2$ where $Y_1$ and
$Y_2$ are Hermitian. Then by the above bound,
\be\label{ineq7}
\| Y \|_{NC} \le \| Y_1 \|_{NC} + \| Y_2  \|_{NC}
\le 2^{2 - 1/p} \, \| Y_1 \|_{CL} + 2^{2 - 1/p} \, \| Y_2 \|_{CL}
\ee
Consider now the
definition of $\| Y_1 + i Y_2 \|_{CL}$.
There are matrices
$Z_0,Z_1,Z_2,Z_3$ such that
\be
\| Y_1 + i Y_2 \|_{CL} = \frac{1}{2}
\, (
\Psi(W) + \Psi(V) )
\ee
where
\be
W &=& Z_0 \ot I_{2} + Z_3 \ot
\sigma_z + (Z_1 + Y_1) \ot \sigma_x
- (Z_2 + Y_2) \ot \sigma_y,
\nonumber \\
V &=& Z_0 \ot I_{2} + Z_3 \ot \sigma_z + Z_1  \ot
\sigma_x
- Z_2 \ot \sigma_y
\ee
Define
\be
{\tilde W}  &=& Z_0 \ot
I_{2} - Z_3 \ot \sigma_z + (Z_1 + Y_1) \ot \sigma_x
+ (Z_2 + Y_2) \ot
\sigma_y, \nonumber \\
{\tilde V} &=& Z_0 \ot I_{2} - Z_3 \ot
\sigma_z + Z_1  \ot \sigma_x
+ Z_2 \ot \sigma_y
\ee
Since $W,V$ are
positive, and since $\tilde W$ and $\tilde V$ are obtained by
conjugation with the unitary $I_{mn} \ot \sigma_x$, it follows that
$\tilde W$ and $\tilde V$
are also positive.
Then  invariance of
$\Psi$ under local unitaries implies that
\be
\Psi(W) = \Psi({\tilde
W}), \quad
\Psi(V) = \Psi({\tilde V})
\ee
By convexity of $\Psi$ it
follows that
\be
\| Y_1 + i Y_2 \|_{CL} \ge \frac{1}{2} \,
(
\Psi((W+{\tilde W})/2) + \Psi((V+{\tilde V})/2) )
\ge \| Y_1 \|_{CL}
\ee
A similar argument shows that $\| Y_1 + i Y_2 \|_{CL} \ge
\| Y_2 \|_{CL}$, hence
from (\ref{ineq7}) it follows that
\be\label{ineq8a}
\| Y \|_{NC} \le 2^{2 - 1/p} \, ( \| Y_1 \|_{CL} +  \| Y_2 \|_{CL})
\le 2^{3 - 1/p} \,  \| Y_1 + i Y_2\|_{CL} = 2^{3 - 1/p} \, \| Y \|_{CL}
\ee

\medskip
\par\noindent \fbox{The case $p > q$}
As in the previous case it is sufficient to consider a Hermitian matrix
$Y = Z - W$ with $Z,W \ge 0$, and to show the analog of (\ref{ineq6a}), that is
\be\label{ineq6b}
\| Y \|_{NC} \le 2^{2 - 1/p} \,\Psi_{p,q} \bmat {Z & 0 \cr 0 & W} \emat
\ee
Once (\ref{ineq6b}) is shown,
the argument leading from (\ref{ineq7}) to (\ref{ineq8a}) can be repeated and
the result then follows.
In order to show (\ref{ineq6b}), we
restrict the infimum on the right side of (\ref{Pis2}) to obtain
\be
\| Y \|_{NC} \le \inf_{C > 0, \, {\rm Tr}\, C =1} \, \| (C^{-1/2r} \ot I_m) Y (C^{-1/2r} \ot I_m) \|_p
\ee
The steps leading from (\ref{first:p<q}) to (\ref{ineq4}) can now be repeated,
leading to the conclusion
\be\label{ineq5}
\| Y \|_{NC}^p \le 2^{2p - 1} \, \inf_{C > 0, \, {\rm Tr}\, C =1} \, \tr
\Big[(C^{-1/2r} \otimes I_m) \bmat {Z & 0 \cr 0 & W} \emat (C^{-1/2r} \otimes I_m)\Big]^p
\ee
Since $\bmat {Z & 0 \cr 0 & W} \emat$ is positive semidefinite, we may
use the Lieb-Thirring bound, as we did for (\ref{Pis2.5}), to conclude that
\be
\| Y \|_{NC}^p \le 2^{2p - 1} \, \Psi_{p,q} \bmat {Z & 0 \cr 0 & W} \emat^p
\ee
 and this completes the proof.

\section{Proof of Theorem \ref{thm2}}\label{sect:proofb}
Next we demonstrate (\ref{thm2:ineq1}) with a family of examples.
There are $2^n$ diagonal $n \times n$ matrices with  $\pm 1$ on the diagonal.
Denote these unitary matrices as $\{U_a\}$, with $U_1 = I_n$.
Then for any $n \times n$ matrix $A$,
\be\label{unit-diag}
\sum_{a=1}^{2^n} 2^{-n} \, U_a \, A \, U_a = A_{diag}
\ee
where $A_{diag}$ is the diagonal matrix obtained by replacing all off-diagonal
entries of $A$ with zero.

Let $| \psi \ket \in \mathbb{C}^d$ be a unit vector. Define
\be
Y_1 = | \psi \ket \bra \psi |, \quad
Y_a = U_a \, Y_1 \, U_a, \quad a=1,\dots,2^n
\ee
Note that $Y_a$ is a pure state for all $a$ and hence $Y_a^p = Y_a$ for all $p$.
Let $\lambda_1,\dots,\lambda_n$ be the diagonal entries of $Y_1$, then define
\be
D =
\sum_{a=1}^{2^n} 2^{-n} \, Y_a = (Y_{1})_{diag} =
\bmat {\lambda_1 & 0 & \dots & 0 \cr 0 & \lambda_2 &  & \cr \vdots && \ddots & \cr 0 &&& \lambda_n} \emat
\ee

\medskip
Let $m = 2^n$ and define the $mn \times mn$ block diagonal matrix
\be
Y = \sum_{a=1}^m Y_k \ot | a \ket \bra a |   =
\bmat {Y_1 & 0 &  & \cdots \cr
0 & Y_2 & 0 & \cdots \cr
\vdots & & \ddots & \cr
0 & & & Y_m} \emat
\ee
Then $Y \ge 0$ and $Y^p=Y$, hence for
$1 \le p \le 2$ and $q \ge 1$
\be\label{CL3}
\Psi_{p,q}(Y)^p =  \bigg( \tr_1 ( \tr_2 Y^p)^{q/p} \bigg)^{p/q} =
\bigg( \tr \Big( \sum_{a=1}^m Y_a \Big)^{q/p} \bigg)^{p/q} =
2^n \, \| D \|_{q/p}
\ee
Furthermore, for $p < q$
\be\label{Pis3}
\| Y \|_{NC}^p & = &
\sup_{C \ge 0, \, {\rm Tr} \, C =1} \, \| (C^{1/2r} \otimes I_m) Y ( C^{1/2r} \otimes I_m) \|_{p}^p \nonumber \\
& = & \sup_{C \ge 0, \, {\rm Tr} \, C =1} \, \sum_{a=1}^m \tr (C^{1/2r}  Y_a  C^{1/2r})^p
\ee
where $r^{-1} = p^{-1} - q^{-1}$.
Note that for any $b \in \{1,\dots,2^n\}$,
\be\label{conj-inv}
\{Y_a\}_{a=1}^{2^n} = \{U_b Y_a U_b\}_{a=1}^{2^n}
\ee
that is conjugation by $U_b$ permutes the set of matrices $Y_a$.
The property (\ref{conj-inv}) implies that for all $C$ and all $b \in \{1,\dots,2^n\}$,
\be\label{inv1}
\sum_{a=1}^m \tr (C^{1/2r}  Y_a  C^{1/2r})^p & = &
\sum_{a=1}^m \tr (C^{1/2r}  U_b Y_a U_b  C^{1/2r})^p \nonumber \\
& = & \sum_{a=1}^m \tr ((U_b C U_b)^{1/2r}  Y_a  (U_b C U_b)^{1/2r})^p
\ee
Furthermore since $r^{-1} \le p/r \le 1$ the map
\be
C \mapsto \tr (C^{1/2r}  Y  C^{1/2r})^p = \tr (Y^{1/2} C^{1/r} Y^{1/2})^p
\ee
is concave on the positive matrices \cite{CL2}. Together with (\ref{inv1}) and (\ref{unit-diag}) this implies that
for all positive $C$,
\be
\sum_{a=1}^m \tr (C^{1/2r}  \,\, Y_a \,\, C^{1/2r})^p \le
\sum_{a=1}^m \tr (C_{diag}^{1/2r}  \,\, Y_a \,\, C_{diag}^{1/2r})^p
\ee
Hence the supremum in (\ref{Pis3}) is achieved on diagonal matrices. Therefore
\be\label{Pis4}
\| Y \|_{NC}^p = \sup_{\{p_1,\dots,p_n \ge 0\}} \, \sum_{a=1}^m  \Big(\sum_{j=1}^n p_j^{1/r}\, \lambda_j \Big)^p
= 2^n \, \sup_{\{p_1,\dots,p_n \ge 0\}} \, \Big(\sum_{j=1}^d p_j^{1/r}\, \lambda_j \Big)^p
\ee
where the $\sup$ runs over positive vectors satisfying  $\sum_j p_j =1$.
Using H\"older's inequality gives
\be
\sum_{j=1}^n p_j^{1/r}\, \lambda_j \le \Big( \sum_{j=1}^d\, \lambda_j^{r'} \Big)^{1/r'}
\ee
where $r'$ is conjugate to $r$. Therefore
\be\label{Pis5}
\| Y \|_{NC}^p \le 2^n \,\Big( \sum_{j=1}^n\, \lambda_j^{r'} \Big)^{p/r'} = 2^n \, \| D \|_{r'}^{p}
\ee
Combining with (\ref{CL3}) gives
\be\label{ratio1}
\frac{\Psi_{p,q}(Y)^p}{\| Y \|_{NC}^p} \ge \frac{\| D \|_{q/p}}{\| D \|_{r'}^{p}}
\ee

Now consider the values
\be
\lambda_j = \frac{c}{j}, \quad j=1,\dots,n
\ee
The normalization condition $\sum \lambda_j =1$ implies that
\be
c \le \frac{1}{\ln n}
\ee
For $t > 1$ define
\be
h(t) = \bigg( \sum_{k=1}^{\infty} k^{-t} \bigg)^{1/t}
\ee
then it follows that
\be
c \le \| D \|_{t} \le c \, h(t)
\ee
Since $p < q$ it follows that $r < \infty$ and $r' > 1$, and
therefore
\be
\frac{\| D \|_{q/p}}{\| D \|_{r'}^{p}} \ge c^{1-p} \, \frac{1}{h(r')^p} \ge (\ln n)^{p-1} \, \frac{1}{h(r')^p}
\ee
Since $p > 1$, $(\ln n)^{p-1}$ diverges as $n \rightarrow \infty$.
Therefore the left side of (\ref{ratio1}) is not uniformly bounded, and we have
a sequence of positive semidefinite matrices $\{Y^{(k)}\}$ such that
\be\label{diverge1}
\frac{\Psi_{p,q}(Y^{(k)})}{\| Y^{(k)} \|_{NC}} \rightarrow \infty
\ee
as $k \rightarrow \infty$. This proves (\ref{thm2:ineq1}).

\medskip
In order to prove (\ref{thm2:ineq2}) we use the result of Lemma \ref{lem:pos-bound5} below,
which shows that for $p=2$ and $q \ge 2$ we can lower bound
$\| Y \|_{CL}$ by $2^{-1/2} \Psi_{p,q}(Y)$. Inserting this bound in (\ref{diverge1})
completes the proof.

\qed

\medskip
\begin{lemma}\label{lem:pos-bound5}
For $2 = p \le q$ and for all $Y \ge 0$,
\be
\| Y \|_{CL} \ge \frac{1}{\sqrt{2}} \, \Psi_{p,q}(Y)
\ee
\end{lemma}

\par\noindent{\em Proof:}
Recall that for a positive Hermitian matrix $Y$,
\be
\| Y \|_{CL} & = & \inf_{A \ge 0} \,\, \Psi_{p,q} \bmat {Y + A
& 0 \cr 0 & A} \emat \nonumber \\
& = & \inf_{A \ge 0} \,\, \bigg[ \tr_1 \Big( \tr_2 (Y+A)^p + \tr_2 A^p \Big)^{q/p} \bigg]^{1/q} \nonumber \\
& =& \inf_{A \ge 0} \,\, \|\tr_2 (Y+A)^p + \tr_2 A^p\|_{q/p}^{1/p}
\ee
Note that for $p=2$,
\be
(Y+A)^2 +  A^2 & = & Y^2 + 2 \, \Big(A + \frac{Y}{2}\Big)^2 - \frac{Y^2}{2} \nonumber \\
& \ge & \frac{Y^2}{2}
\ee
The partial trace preserves positivity, and since $q \ge 2$ the $q/2$-norm is monotonic, hence
\be
\|\tr_2 (Y+A)^2 + \tr_2 A^2\|_{q/2} \ge
 \frac{1}{2} \,  \| \tr_2 Y^2 \|_{q/2}
= \frac{1}{2} \,  \Psi_{2,q}(Y)^2
\ee
as claimed.

\bigskip

  \noindent{\bf Acknowledgments:}
This work  was  supported in part by the
National Science Foundation under grant DMS-0400426. The authors thank E. Carlen
for helpful suggestions and discussions.

\pagebreak

{~~}

\pagebreak
\appendix

\section*{Appendix A: Proof of Lemma \ref{lem1}}
Because the norm is given by the different expressions (\ref{Pis1}) and (\ref{Pis2}) depending
on the relative sizes of $p$ and $q$, the proofs are given separately for the two cases.

\bigskip

\par\noindent \fbox{{\bf (a)} : $p \le q$}
For any matrices $A,B \in {\cal M}_n$ and $Y,W \in {\cal M}_{nm}$,
\bee
\| (A \ot I_m) (Y + W) (B \ot I_m) \|_p & \le & \| (A \ot I_m) Y (B \ot I_m) \|_p \\
&& \hskip0.4in + \hskip0.2in \| (A \ot I_m) W (B \ot I_m) \|_p \\
& \le & \| A\|_{2r} \, \| B \|_{2r} \, \Big( \| Y \|_{NC} + \| W \|_{NC} \Big)
\eee
Dividing both sides by $\| A\|_{2r} \, \| B \|_{2r}$ and taking the $\sup$ over $A,B$ gives the bound.

\medskip
\medskip
\par\noindent \fbox{{\bf (a)} : $q \le p$}
Given $\epsilon > 0$ there are matrices $A_i, B_i, Z_i$ such that
\bee
Y = \sum_i (A_i \ot I_m) Z_i (B_i \ot I_m), \qquad
\sum_i \|A_i\|_{2r} \, \|B_i\|_{2r} \, \|Z_i\|_{p} < \| Y \|_{NC} + \epsilon
\eee
Similarly there are matrices $C_j,D_j,X_j$ such that
\bee
W = \sum_j (C_j \ot I_m) X_j (D_j \ot I_m), \qquad
\sum_j \|C_j\|_{2r} \, \|D_j\|_{2r} \, \|X_j\|_{p} < \| W \|_{NC} + \epsilon
\eee
 Therefore
\bee
\| Y + W \|_{NC} & \le & \sum_i \|A_i\|_{2r} \, \|B_i\|_{2r} \, \|Z_i\|_{p} +
\sum_j \|C_j\|_{2r} \, \|D_j\|_{2r} \, \|X_j\|_{p}  \\
&<&
\| Y \|_{NC} + \| W \|_{NC} + 2 \epsilon
\eee
Since this holds for all $\epsilon > 0$ the result follows.

\medskip
\medskip
\par\noindent \fbox{{\bf (b)}}
Without loss of generality assume that $p \le q$, then it follows that $q' \le p'$. Consider any decomposition
of $W$:
\be\label{ineq:b1}
W = \sum_i \, (A_i \ot I_m) Z_i (B_i \ot I_m)
\ee
Hence
\bee
| \tr (Y W) | & = & | \sum_i \, \tr [Y (A_i \ot I_m) Z_i (B_i \ot I_m) ] | \\
& = & | \sum_i \, \tr [ Z_i (B_i \ot I_m) Y (A_i \ot I_m) ] | \\
& \le & \sum_i \| Z_i \|_{p'} \, \| (B_i \ot I_m) Y (A_i \ot I_m) \|_{p} \\
& \le & \sum_i \, \| Z_i \|_{p'} \, \| B_i \|_{2 r} \, \| A_i \|_{2 r} \, \| Y \|_{NC:p,q}
\eee
Noting that $r^{-1} = p^{-1} - q^{-1} = (q')^{-1} - (p')^{-1}$  we may take the $\inf$ over
$A_i,B_i,Z_i$ satisfying (\ref{ineq:b1}) to conclude that
\bee
| \tr (Y W) | \le \| Y \|_{NC:p,q} \, \| W \|_{NC:p',q'}
\eee

\medskip
\medskip
\par\noindent \fbox{{\bf (c)} : $p \le q$}
From the definition (\ref{Pis1}) we obtain
\bee
\| Y \|_{NC:p,q} & = & \sup_{A,B} \bigg\{ \| (A \ot I_m) Y (B \ot I_m) \|_{p} \,:\, \| A \|_{2r} \, \| B \|_{2r} \le 1 \bigg\} \\
& = & \sup_{A,B,Z} \, \bigg\{ | \tr [(A \ot I_m) Y (B \ot I_m) Z] | \,:\, \| A \|_{2r} \, \| B \|_{2r} \le 1, \, \| Z \|_{p'} \le 1 \bigg\} \\
& \le & \sup_{A,B,Z} \, \bigg\{ | \tr [ Y (B \ot I_m) Z (A \ot I_m)] | \,:\, \| A \|_{2r} \, \| B \|_{2r} \, \| Z \|_{p'} \le 1 \bigg\} \\
& = & \sup_{W,A,B,Z} \, \bigg\{ | \tr ( Y W) | \,:\, W = (B \ot I_m) Z (A \ot I_m), \\
&& \hskip1.5in \,\, \| A \|_{2r} \, \| B \|_{2r} \, \| Z \|_{p'} \le 1 \bigg\} \\
& \le & \sup_{W} \, \bigg\{ | \tr (Y W ) | \,:\, \| W \|_{NC:p',q'} \le 1 \bigg\}
\eee
where in the last inequality we used again $r^{-1} = (q')^{-1} - (p')^{-1}$. Using H\"older's inequality,
the condition  $\| W \|_{NC:p',q'} \le 1$ implies
\bee
| \tr (Y W ) | \le  \| Y \|_{NC:p,q} \,\| W \|_{NC:p',q'} \le  \| Y \|_{NC:p,q}
\eee
Combining these inequalities we deduce that equality must hold, and hence the result follows.

\medskip
\medskip
\par\noindent \fbox{{\bf (c)} : $q \le p$}
From H\"older's inequality we know that
\bee
\sup_{W} \, \bigg\{ |
\tr (Y W ) | \,:\, \| W \|_{NC:p',q'} \le 1 \bigg\} \le  \| Y
\|_{NC:p,q}
\eee
so it is sufficient to show that there is a
matrix $W$ such that
\be\label{c:q<p1} |
\tr (Y W ) | = \| Y
\|_{NC:p,q} \,\, \| W \|_{NC:p',q'}
\ee
Consider the space
$X={\cal M}_{nm}$ equipped with the norm (\ref{Pis2}). Every
linear functional on $X$ may be written as $f(\cdot) = \tr (\cdot
\, W)$ for some matrix $W$, therefore the norm of $f$ is
\bee
\| f
\| = \sup_{Z}  \{ | \tr (W Z) | \,:\, \| Z \|_{NC:p,q} \le 1 \} =
\| W \|_{NC:p',q'}
\eee
where we used the previous duality result. So
it is sufficient to show that for every $Y \in X$ there is a
linear functional $f$ with $f(Y)  = \| f \| \, \| Y \|_{NC:p,q}$,
as this will imply (\ref{c:q<p1}). But the existence of such a functional is a well-known
corollary of the Hahn-Banach Theorem (see for example \cite{RS1}).

\medskip
\medskip
\par\noindent \fbox{{\bf (d)} : $p \le q$}
First note that for any matrices $A,B$
\bee
\| (A \ot I_m) (Y_1 \ot Y_2) (B \ot I_m) \|_p =
\| A Y_1 B \|_p \, \| Y_2 \|_p
\eee
Furthermore $p^{-1} = r^{-1} + q^{-1}$ hence by H\"older's inequality
\bee
\| A Y_1 B \|_p \, \| Y_2 \|_p \le \| A \|_{2r} \, \| Y_1 \|_{q} \, \| B \|_{2r} \, \| Y_2 \|_p
\eee
Therefore
\be\label{ineq:b,p<q}
\frac{\| (A \ot I_m) (Y_1 \ot Y_2) (B \ot I_m) \|_p}{\| A \|_{2r} \,  \| B \|_{2r}}
\le \| Y_1 \|_{q} \, \| Y_2 \|_p
\ee
and taking the $\sup$ over $A,B$ gives
\bee
\| Y_1 \ot Y_2 \|_{NC} \le \| Y_1 \|_{q} \, \| Y_2 \|_p
\eee
It remains to show that equality can be achieved in (\ref{ineq:b,p<q}) by a suitable
choice of $A,B$. Without loss of generality we can assume that $Y_1 \ge 0$, in which case we define
\bee
A = B = Y_1^{(q-p)/2p}
\eee
Then $A Y_1 B = Y_1^{q/p}$, and $\| A \|_{2r} = \| B \|_{2r} = \| Y_1 \|_{q}^{q/2r}$, which yields
\bee
\frac{\| (A \ot I_m) (Y_1 \ot Y_2) (B \ot I_m) \|_p}{\| A \|_{2r} \,  \| B \|_{2r}} =
\| Y_1^{q/p} \|_{p} \, \| Y_1 \|_{q}^{- q/r} \,  \| Y_2 \|_p
=   \| Y_1 \|_{q} \,  \| Y_2 \|_p
\eee

\medskip
\medskip
\par\noindent \fbox{{\bf (d)} : $q \le p$}
There are matrices $W_1, W_2$ satisfying
\bee
\| Y_1 \|_{q} & = &  \tr Y_1 \, W_1, \quad \| W_1 \|_{q'} = 1 \\
\| Y_2 \|_p & = & \tr Y_2 \, W_2, \quad \| W_2 \|_{p'} = 1
\eee
where $q', p'$ are the conjugate values for $q,p$ respectively.
Note that $p' \le q'$.
Consider any set of matrices $A_i,B_i, Z_i$ for which
\be\label{ineq:b,q<p1}
Y_1 \ot Y_2 = \sum_i (A_i \ot I_m) Z_i (B_i \ot I_m)
\ee
then it follows that
\bee
\| Y_1 \|_{q} \,  \| Y_2 \|_p & = & \tr (Y_1 \ot Y_2) (W_1 \ot W_2) \\
& = & \sum_i \tr (A_i \ot I_m) Z_i (B_i \ot I_m) (W_1 \ot W_2) \\
& = & \sum_i \tr Z_i \Big( (B_i W_1 A_i) \ot W_2 \Big) \\
& \le & \sum_i \| Z_i \|_{p} \, \| B_i W_1 A_i \|_{p'} \, \| W_2 \|_{p'} \\
& \le & \sum_i \| Z_i \|_{p} \, \| B_i \|_{2 r'} \, \| W_1 \|_{q'} \, \| A_i \|_{2 r'} \, \| W_2 \|_{p'} \label{ineq:b,q<p2}
\eee
where in the last line the value $r'$ is given by
\bee
\frac{1}{r'} = \frac{1}{p'} - \frac{1}{q'} = 1 - \frac{1}{p} - \Big( 1 - \frac{1}{q} \Big) = \frac{1}{r}
\eee
Hence $r' = r$ and therefore (\ref{ineq:b,q<p2}) gives
\bee
\| Y_1 \|_{q} \,  \| Y_2 \|_p \le  \sum_i \| Z_i \|_{p} \, \| B_i \|_{2 r} \,  \| A_i \|_{2 r}
\eee
Since this holds for every decomposition  (\ref{ineq:b,q<p1}) it follows that
\bee
\| Y_1 \|_{q} \,  \| Y_2 \|_p \le \| Y_1 \ot Y_2 \|_{NC}
\eee
To complete the proof, it is sufficient to find matrices $A,B, Z$ such that $Y_1 = A Z B$ and
\bee
\| Y_1 \|_{q} = \| A \|_{2 r} \, \| Z \|_{p} \, \| B \|_{2 r}
\eee
Without loss of generality we may again assume that $Y_1 \ge 0$, then taking
\bee
A = B = Y_1^{(p-q)/2p}, \quad Z = Y_1^{q/p}
\eee
it follows that $A Z B = Y_1$ and $\| A \|_{2r}^2 = \| Y_1 \|_q^{q/r}$, and hence
\bee
\| A \|_{2 r} \, \| Z \|_{p} \, \| B \|_{2 r} = \| Y_1 \|_q^{q/r} \, \| Z \|_{p} = \| Y_1 \|_{q}^{q/r} \, \| Y_1 \|_{q}^{q/p}
= \| Y_1 \|_q
\eee

\medskip
\medskip
\par\noindent \fbox{{\bf (e)} : $p \le q$}
We must show that the expression (\ref{Pis1}) reduces to $\| Y \|_p$ when $p=q$.
In this case $r = \infty$, and so
\bee
\| Y \|_{NC} = \sup_{A,B} \frac{\| (A \ot I_m) Y (B \ot I_m) \|_p}{\| A \|_{\infty} \, \| B \|_{\infty}}
\eee
Note that $\| A \ot I_m \|_{\infty} = \| A \|_{\infty}$ and similarly for $B$, therefore by H\"older's inequality
\bee
\| (A \ot I_m) Y (B \ot I_m) \|_p \le \| A \ot I_m \|_{\infty} \, \| B \ot I_m \|_{\infty} \, \| Y \|_p
= \| A  \|_{\infty} \, \| B  \|_{\infty} \, \| Y \|_p
\eee
which leads to
\be\label{ineq:c,p<q}
\| Y \|_{NC} \le \| Y \|_p
\ee
Taking $A =  B = I_n$ shows that equality is achieved in (\ref{ineq:c,p<q}).

\medskip
\medskip
\par\noindent \fbox{{\bf (e)} : $q \le p$}
Now we show that (\ref{Pis2}) also reduces to $\| Y \|_p$ when $p=q$.
Consider any decomposition $Y = \sum_i (A_i \ot I_m) Z_i (B_i \ot I_m)$ then
\bee
\| Y \|_p & \le & \sum_i \| (A_i \ot I_m) Z_i (B_i \ot I_m) \|_p \\
& \le & \sum_i \| A_i \ot I_m \|_{\infty} \, \| Z_i \|_p \, \| B_i \ot I_m \|_{\infty} \\
& = & \sum_i \| A_i  \|_{\infty} \, \| Z_i \|_p \, \| B_i  \|_{\infty}
\eee
 and hence
\bee
\| Y \|_p \le \| Y \|_{NC}
\eee
Taking $A_1 = B_1 = I_n$ and $Z_1 = Y$ shows that equality is achieved, and hence the result follows.

\medskip
\medskip
\par\noindent \fbox{{\bf (f)} : $p \le q$}
Let $Y$ be diagonal in the product basis with entries $\{ y_{ij} \}$. Let $A \in {\cal M}_{n}$ be the diagonal matrix with entries
\be\label{def:A:f,p<q}
a_i = w_i^{q/2pr} \,\, \bigg( \sum_i w_i^{q/p} \bigg)^{\frac{-1}{2r}}
\ee
where $w_i = \sum_j | y_{ij} |^p$. Then $\| A \|_{2r} = (\sum_i | a_i |^{2r})^{1/{2r}} = 1$ and
\be
\tr | (A \ot I_m) Y (A \ot I_m)|^p &=& \sum_{i} | a_i |^{2p} \, w_i = \bigg( \sum_i w_i^{q/p} \bigg)^{p/q}\\
&=& \bigg( \tr_1 ( \tr_2 |Y|^p)^{q/p} \bigg)^{p/q}
\ee
Hence we deduce that
\be\label{f:p<q3}
\| Y \|_{NC} \ge \bigg( \tr_1 ( \tr_2 |Y|^p)^{q/p} \bigg)^{1/q}
\ee

Furthermore, write $Y = | Y| V$ where $V$ is a diagonal unitary matrix, then for any $A,B$
\bee
\| (A \ot I_m) \hskip-0.1in &Y&  \hskip-0.1in (B \ot I_m) \|_{p}   \\
& = & \| (A \ot I_m) |Y|^{1/2} |Y|^{1/2} V (B \ot I_m) \|_{p} \\
& \le & \| (A \ot I_m) |Y| (A^* \ot I_m) \|_{p}^{1/2} \,\, \| (B^* \ot I_m) V^* |Y| V (B \ot I_m) \|_{p}^{1/2} \\
& = & \| (A \ot I_m) |Y| (A^* \ot I_m) \|_{p}^{1/2} \,\, \| (B^* \ot I_m) |Y| (B \ot I_m) \|_{p}^{1/2}
\eee
since $V$ is diagonal and therefore commutes with $|Y|$. Therefore we deduce that
\be
\| Y \|_{NC} & \le & \sup_{A} \, \frac{\| (A \ot I_m) |Y| (A^* \ot I_m) \|_{p}}{\| A \|_{2r}^2} \\
& = &
\sup_{A \ge 0} \, \frac{\| (A \ot I_m) |Y| (A \ot I_m) \|_{p}}{\| A \|_{2r}^2} \\
& = & \sup_{\sigma} \, \bigg( \tr \bigg( (\sigma^{1/2r} \ot I_m) |Y| (\sigma^{1/2r} \ot I_m) \bigg)^p \bigg)^{1/p}\\
& = & \sup_{\sigma} \, \bigg( \tr \bigg( |Y|^{1/2} \, (\sigma^{1/r} \ot I_m) |Y|^{1/2} \bigg)^p \bigg)^{1/p} \label{f:p<q2}
\ee
where the final $\sup$ runs over positive semidefinite matrices with trace one.
Define
\be
F(\sigma) = \tr \bigg( |Y|^{1/2} \, (\sigma^{1/r} \ot I_m) |Y|^{1/2} \bigg)^p
\ee
Since $r \ge p$ the
results of Theorem 1.1 in \cite{CL2} imply that $F$ is concave.
Furthermore there are $2^n$ diagonal matrices $\{ U_j \}$ with $\pm 1$ on the diagonals such that
\be
(\sigma)_{diag} = 2^{-n} \, \sum_j U_j \, \sigma \, U_j^*
\ee
Since $|Y|$ is also diagonal this implies that
\bee
F((\sigma)_{diag}) & = & F( 2^{-n} \, \sum_j U_j \, \sigma \, U_j^* ) \\
& \ge & 2^{-n} \, \sum_j \, F(U_j \, \sigma \, U_j^* ) \\
& = & 2^{-n} \, \sum_j \, F( \sigma ) \\
& = & F(\sigma)
\eee
Hence the $\sup$ in (\ref{f:p<q2}) is achieved on the diagonal matrices.
Therefore
\be
\| Y \|_{NC} & \le & \sup_{c_i \ge 0} \, \bigg( \sum_{i,j} |c_i y_{ij} |^p \bigg)^{1/p}
= \sup_{c_i \ge 0} \, \bigg( \sum_{i} | c_i |^p \, w_i \bigg)^{1/p}
\ee
where the $\{ c_i \}$ satisfy $\sum_i | c_i |^r = 1$.
Using H\"older's inequality with $s = r/p, s' = q/p$ gives
\bee
\sum_{i} | c_i |^p \, w_i & \le & \bigg( \sum_i | c_i |^r \bigg)^{1/s} \,\, \bigg( \sum_i w_i^{s'} \bigg)^{1/s'} \\
& \le & \bigg( \sum_i w_i^{q/p} \bigg)^{p/q} \\
& = & \bigg( \tr_1 ( \tr_2 |Y|^p)^{q/p} \bigg)^{p/q}
\eee
and therefore
\be\label{f:p<q4}
\| Y \|_{NC} & \le & \bigg( \tr_1 ( \tr_2 |Y|^p)^{q/p} \bigg)^{1/q}
\ee
Combining (\ref{f:p<q3}) and (\ref{f:p<q4}) gives the result.

\medskip
\medskip
\par\noindent \fbox{{\bf (f)} : $q \le p$}
By using the same diagonal matrix $A$ defined for the $p \le q$ case in (\ref{def:A:f,p<q}), we can write
$Y = (A \ot I_m) Z (A \ot I_m)$. Then we have $\| A \|_{2r}=1$ and
\be
\| Z \|^p_p &=& \sum_{i,j}|z_{ij}|^p = \sum_i \frac{w_i}{|a_i|^{2p}}= \bigg( \sum_i w_i^{q/p} \bigg)^{p/q}\\
&=& \bigg( \tr_1 ( \tr_2 |Y|^p)^{q/p} \bigg)^{p/q}
\ee
Therefore
\be\label{f:q<p1}
\| Y \|_{NC} \le \bigg( \tr_1 ( \tr_2 |Y|^p)^{q/p} \bigg)^{1/q}
\ee
Furthermore, by classical duality there is a diagonal matrix $W$ satisfying
\be
\bigg( \tr_1 ( \tr_2 |W|^{p'})^{q'/p'} \bigg)^{1/q'} =1, \qquad
\bigg( \tr_1 ( \tr_2 |Y|^p)^{q/p} \bigg)^{1/q} = | \tr ( YW ) |
\ee
Consider any decomposition $Y = \sum_i (A_i \ot I_m) Z_i (B_i \ot I_m)$, then
\bee
| \tr (YW) | & = & | \sum_i \, \tr [ Z_i (B_i \ot I_m) W (A_i \ot I_m) ] | \\
& \le &  \sum_i \, \| Z_i \|_p \, \| (B_i \ot I_m) W (A_i \ot I_m) \|_{p'} \\
& \le & \sum_i \, \| Z_i \|_p \, \| B_i \|_{2r'} \, \| A_i \|_{2r'} \, \| W \|_{NC:p',q'} \\
& = & \sum_i \, \| Z_i \|_p \, \| B_i \|_{2r'} \, \| A_i \|_{2r'}
\eee
where in the last line we used the previous result that for diagonal matrices
\be
\| W \|_{NC:p',q'} = \bigg( \tr_1 ( \tr_2 |W|^{p'})^{q'/p'} \bigg)^{1/q'}
\ee
since $p' \le q'$. Also $r' = r$ and so
\be
| \tr (YW) | \le \sum_i \, \| Z_i \|_p \, \| B_i \|_{2r} \, \| A_i \|_{2r}
\ee
Since this holds for every decomposition we deduce that
\be\label{f:q<p2}
\bigg( \tr_1 ( \tr_2 |Y|^p)^{q/p} \bigg)^{1/q} \le \| Y \|_{NC}
\ee
Together (\ref{f:q<p1}) and (\ref{f:q<p2}) imply equality.

\section*{Appendix B: Proof of Lemma \ref{lem:non-mon}}
\par\noindent{\em Proof:}
Let $Y \ge 0$ be a positive semidefinite matrix. Recall Lemma \ref{lem2}. Setting 
$A=0$ on the right side of (\ref{eqn:lem1})
shows that
\be\label{eqn14}
\| Y \|_{CL} \le \Psi_{p,q} \bmat {Y & 0 \cr 0 & 0} \emat =
\Psi_{p,q}(Y)
\ee
Furthermore, for all $A \ge 0$,
\be\label{eqn14a}
\Psi_{p,q} \bmat {Y + A & 0 \cr 0 & A} \emat^q =
\tr_1 \bigg(\tr_2 (Y+A)^p + \tr_2 A^p\bigg)^{q/p}
\ee
The Schatten norm $A \mapsto \| A \|_t$ is monotone for all $t \ge 1$.
Thus the condition $q \ge p$ implies from (\ref{eqn14a}) that
\be
\Psi_{p,q} \bmat {Y + A & 0 \cr 0 & A} \emat^q
\ge \tr_1 \bigg(\tr_2 A^p\bigg)^{q/p} 
= \Psi_{p,q}(A)^q
\ee
Hence the infimum on the right side of (\ref{eqn:lem1}) may be 
restricted to positive matrices $A$
satisfying $\Psi_{p,q}(A) \le \Psi_{p,q}(Y)$. This is a compact set, 
and the function
$\Psi_{p,q}$ is continuous, hence the infimum is achieved on a 
positive semidefinite matrix $B$. 
Furthermore if $C,D \ge 0$ are nonzero positive semidefinite matrices
and $t \ge 1$ then
$\| C \|_t < \| C + D \|_t$ (this can be seen by reducing to the case where $D$
has rank one, and then using the interlacing condition for the eigenvalues of
$C$ and $C+D$). If we assume $B \neq 0$ then this implies
\be\label{eqn15}
\| Y \|_{CL} & = & \Psi_{p,q} \bmat {Y+B & 0 \cr 0 & B} \emat \nonumber \\
& > & \Psi_{p,q} \bmat {Y+B & 0 \cr 0 & 0} \emat \nonumber \\
& = & \Psi_{p,q}(Y+B)  \nonumber \\
& \ge & \| Y+B \|_{CL}
\ee
where in the first inequality we used strict monotonicity of the map
$A \mapsto \| A \|_{q/p}$, and the second inequality follows from
(\ref{eqn14}).
Since $Y \le Y+B$ this shows that $\| \cdot \|_{CL}$ is not 
monotone unless $B=0$ for every $Y \ge 0$.

But if $B=0$ for all $Y \ge 0$ it would follow from (\ref{eqn15}) that $\| Y \|_{CL} = 
\Psi_{p,q}(Y)$
for all $Y \ge 0$. However it is easy to construct examples which show that
the function $Y \mapsto \Psi_{p,q}(Y)$ is not monotone (see below for one such
construction), and therefore also in this case we conclude that 
$\| \cdot \|_{CL}$ is not 
monotone.

\medskip
We now present a numerical example showing that $\Psi_{p,q}$ is not monotone.
Let $w = 3 - \sqrt{10}$ and define the positive matrices
\be
W = \bmat {w^2 & w&0 & 0 \cr w & 1&0 & 0\cr 0&0&0&0 \cr 0&0&0&0} \emat, \quad
Y = \bmat {1&0&0&0 \cr0&0&0&0 \cr 0&0&1/2&1/2 \cr 0&0&1/2&/12} \emat
\ee
Define the function
\be
g(t) = \Psi_{p,q}(Y + t W)
\ee
Then $g'(t) < 0$ for small values of $t$, and a range of values of $p,q$.

\end{document}